\def\editmode{0}

\def\bibfilenames{bibman_refs,other_refs}

\documentclass[twocolumn,twoside]{IEEEtran}

\usepackage{savesym} 
\savesymbol{labelindent}
\usepackage{enumitem} 

\if\editmode1  
\usepackage[backend=bibtex,style=alphabetic,sorting=debug,maxbibnames=99]{biblatex}
\DeclareFieldFormat{labelalpha}{\thefield{entrykey}}
\DeclareFieldFormat{extraalpha}{}
\bibliography{\bibfilenames}
\newcommand{\cmt}[1]{\noindent\textcolor{lightgreen}{\underline{[#1]}}} 
\newcommand{\hc}[1]{\textcolor{blue}{#1}} 

\setlistdepth{9}
\newlist{mydeepitemize}{enumerate}{9}
\setlist[mydeepitemize,1]{label=$\bullet$}
\setlist[mydeepitemize,2]{label=$\diamond$}
\setlist[mydeepitemize,3]{label=$\rightarrow$}
\setlist[mydeepitemize,4]{label=$\circ$}
\setlist[mydeepitemize,5]{label=$-$}
\setlist[mydeepitemize,6]{label=$\square$}
\setlist[mydeepitemize,7]{label=$\star$}
\setlist[mydeepitemize,8]{label=$\checkmark$}
\setlist[mydeepitemize,9]{label=$\Delta$}
\newenvironment{myitemize}{\begin{mydeepitemize}}{\end{mydeepitemize}}

\else
\usepackage{cite}
\newcommand{\cmt}[1]{} 
\newcommand{\hc}[1]{\textcolor{black}{#1}} 

\fi

\newcommand{\printmybibliography}{
\if\editmode1 
\printbibliography
\else
\bibliography{\bibfilenames}
\fi
}

\usepackage{fixltx2e}

\usepackage{graphicx}

\usepackage{subcaption}              

\usepackage[utf8]{inputenc}

\usepackage{amsfonts}
\usepackage{amsmath}
\usepackage{mathtools}

\usepackage{amssymb}

\usepackage{bm}

\usepackage{color,verbatim}
\usepackage{multirow}
\usepackage{accents}

\usepackage{theoremref}

\usepackage{url}

\newcounter{rulecounter}
\newcommand{\resetrule}{ \setcounter{rulecounter}{0}}
\resetrule

\newtheorem{myauxproblem}{Problem}

\newtheorem{myauxoptionalproblem}{Optional Problem}

\newsavebox{\selvestebox}
\newenvironment{colbox}[1]
  {\newcommand\colboxcolor{#1}%
   \begin{lrbox}{\selvestebox}%
   \begin{minipage}{\dimexpr\columnwidth-2\fboxsep\relax}}
  {\end{minipage}\end{lrbox}%
   \begin{center}
   \colorbox{\colboxcolor}{\usebox{\selvestebox}}
   \end{center}}

\definecolor{orange}{rgb}{1,0.8,0}
\definecolor{gray}{rgb}{.9,0.9,0.9}
\definecolor{darkgray}{rgb}{.3,0.3,0.3}
\definecolor{darkblue}{rgb}{.1,0.0,0.3}
\definecolor{lightblue}{rgb}{0.7,0.7,1}
\definecolor{lightred}{rgb}{1,0.7,.7}
\definecolor{purple}{RGB}{204,153,255}
\definecolor{lightgray}{rgb}{.95,0.95,0.95}
\definecolor{lightgreen}{rgb}{0.3,0.5,0.3}
\definecolor{darkgreen}{rgb}{0.05,0.3,0.05}

\newcommand{\brackets}[1]{\left\{#1\right\}}

\newcommand{\rfield}{\mathbb{R}}
\newcommand{\vect}{\mathop{\rm vec}}

\newcommand{\diag}[1]{\mathop{\rm diag}\brackets{#1}}

\newcommand{\rank}{\mathop{\rm rank}}

\newcommand{\range}[1]{\mathcal{R}\brackets{#1}}

\newcommand{\transpose}{^\top}
 
\newtheorem{myproposition}{Proposition}
\newtheorem{myremark}{Remark}
\newtheorem{myproblemstatement}{Problem Statement}
\newtheorem{mylemma}{Lemma}
\newtheorem{mytheorem}{Theorem}
\newtheorem{mydefinition}{Definition}
\newtheorem{mycorollary}{Corollary}

\renewcommand{\range}[1]{\mathcal{R}\{#1\}}

\renewcommand{\vect}{{\rm vec}}

\renewcommand{\transpose}{^\top}

\newcommand{\uparmat}{{\hc{\bm U}}_{\parallel}}

\newcommand{\upermat}{{\hc{\bm U}}_{\perp}}

\usepackage{xr,xcite}
\externalcitedocument[P-]{../paper/draft}

\makeatletter
\renewcommand\@bibitem[1]{\item\if@filesw \immediate\write\@auxout
    {\string\bibcite{#1}{\the\value{\@listctr}}}\fi\ignorespaces}
\def\@biblabel#1{[#1]}
\makeatother

\newcommand{\itembf}[1]{\item\textbf{#1}}
\newcommand{\changed}[1]{{\color{black}#1}}

\begin{document}

\begin{center}
{\Large \textbf{Comments on ``Design of Asymmetric Shift Operators for
Efficient Decentralized Subspace Projection''}}\\[.3cm] Daniel
Romero,\footnote{The author is with the Dept. of Information and
Communication Technology, University of Agder, Grimstad, 4879, Norway
(e-mail: daniel.romero@uia.no).  } \emph{member, IEEE}\\[.3cm]
\end{center}

\begin{abstract}

  This correspondence disproves the main results in the paper ``Design
  of Asymmetric Shift Operators for Efficient Decentralized Subspace
  Projection'' by S. Mollaebrahim and
  B. Beferull-Lozano. Counterexamples and counterproofs are provided
  when applicable. \changed{A correction is suggested for some of the
    flaws. However, it does not seem possible to amend most of the
    flaws since the overall approach based on a Schur decomposition of
    the shift matrix does not appear to be helpful to solve the
    desired problem.  }

\end{abstract}

\section*{Introduction}

This correspondence comments on \cite{siavash2021asymmetric}, where
the goal is to design shift matrices for graph filters whose output is
the result of projecting their input onto a given subspace. The paper
relies on a Schur decomposition
$\bm S = \bm W (\bm D + \bm Q) \bm W\transpose$ of the sought shift
matrix $\bm S$. The goal is therefore the same as in the earlier paper
\cite{romero2021fast} with the exception that $\bm S$ is required to
be symmetric in \cite{romero2021fast}. The reason for such a
requirement in~\cite{romero2021fast} was \changed{to obtain suitable
  necessary and sufficient conditions to characterize the set of
  feasible shift matrices}. The work in \cite{siavash2021asymmetric}
claims to solve the problem when that assumption is lifted. The
present correspondence shows that this is not the case.

\emph{Notation:} The notation in \cite{siavash2021asymmetric} is
adopted here. Besides, $\bm A^0$ denotes the identity matrix, where
$\bm A$ is an arbitrary square matrix.  Equation numbers of the form
($N$) refer to \cite{siavash2021asymmetric}, whereas
(C$N$) refers to the present correspondence.


\section*{Technical Content}
\label{sec:correctness}

\textbf{Problem formulation.} Early in Sec.~IV, the problem is formulated as
finding \emph{asymmetric} matrices that are ``polynomial'' and
``topological''. Although not mentioned
in \cite{siavash2021asymmetric}, these definitions were introduced
by \cite[Sec.~III]{romero2021fast}. The aforementioned formulation
in \cite{siavash2021asymmetric} excludes symmetric matrices, which
conflicts with the optimization approach adopted later. The reason is
that symmetric matrices are included in the closure of the set of
asymmetric matrices and the infimum of a continuous function on a set
equals the infimum of that function on the closure of that set. This first issue
could be amended by replacing ``asymmetric'' with ``possibly
asymmetric''.

\textbf{Theorem 1, Part 1.} Although not mentioned, this result aims
at extending \cite[Th. 1,
C2]{romero2021fast}. In \cite{siavash2021asymmetric}, Theorem 1,
Part 1,
states that ``\emph{A necessary condition to have a valid polynomial
graph shift matrix $\bm S$ is that $\bm D_1$ and $\bm D_2$ do not
share any common value.}'', which is obviously incorrectly phrased as
these matrices are diagonal and therefore they necessarily share the
entry 0, unless they are $1\times 1$. By looking at the proof, it is clear that what
this result aims to say is that the aforementioned two matrices ``do
not share any common eigenvalue'' or, equivalently, ``do not share any
common diagonal value''.

But even such a (corrected) statement would be false, as proved by
the following counterexample. Let
\begin{align}
\label{eq:counterexample1}
\uparmat &= [\bm u_1, \bm u_2],~ \upermat=[\bm u_3, \bm u_4], ~\bm
U=[\uparmat, \upermat],\\ \bm D_1&=\bm
D_2=\diag{[1,0]},~\bm Q=\bm 0,~ \bm W =
[\bm u_1, \bm u_3, \bm u_2, \bm u_4], \nonumber
\end{align}
where $\{\bm u_i\}_i\subset \rfield^4$ are arbitrary orthonormal vectors.  This
clearly satisfies that $\bm S=\bm W(\bm D+\bm Q)\bm
W\transpose=\uparmat\uparmat\transpose$ and therefore it is a
``polynomial shift'' according to Sec.~IV, yet it violates
the aforementioned ``necessary'' condition. One can also obtain an
asymmetric $\bm S$ violating this condition; set e.g. $(\bm Q)_{2,4}$
equal to 1, which yields $\bm S^2=\uparmat\uparmat\transpose$.

Since the result is incorrect, the proof is necessarily incorrect. The
issue can be found where it introduces the  assumption that there exist $ \bm
c, \bm D,\text{ and }\bm Q$ such that
\begin{align}
\sum_{l=0}^{L-1} c_l (\bm
D+\bm Q)^l = 
\left[
\begin{array}{cc}
\bm I_r & \bm 0\\
\bm 0, & \bm 0
\end{array}
\right].
\end{align}
This assumption is not a hypothesis of the theorem and \emph{does} entail
loss of generality. Therefore, this step is logically flawed. It is
clear that one cannot introduce arbitrary assumptions in a
proof. Otherwise one could prove any statement $s$ and its negation
$\bar s$, which would yield a contradiction.

This flaw propagates throughout the paper and invalidates the proposed
optimization problems; cf. the discussions below\-.

Finally, the proof adapts the arguments in \cite[Explanation around
(12)]{romero2021fast} and is also related to those
in \cite{sandryhaila2014finitetime}
and \cite{segarra2017graphfilterdesign}. However, this is presented
without reference to these works.

\changed{Theorem 1, Part 1, tries to extend the necessary (and also
  sufficient) condition in \cite[Theorem 1, C2]{romero2021fast}, to
  asymmetric shift matrices. In \cite{romero2021fast}, this is possible
  because assuming that $\bm S$ is symmetric allows one to use its
  eigendecomposition $\bm S= \bm W\bm \Lambda \bm W\transpose$. The
  columns of $\bm W$ are eigenvectors of $\sum_l c_l \bm S^l$ for all
  $\{c_l\}_l$, which enables one to group the eigenvalues of $\bm S$
  into two groups: those associated with eigenvectors in
  $\range{\uparmat}$ and those associated with eigenvectors in
  $\range{\upermat}$. The necessary condition is that the same
  eigenvalue cannot be in both groups.  In contrast, the columns of
  $\bm W$ in the Schur decomposition
  $\bm S = \bm W (\bm D + \bm Q)\bm W\transpose$ adopted in
  \cite{siavash2021asymmetric} are not generally eigenvectors of
  $\sum_l c_l \bm S^l$ and, therefore, no such a grouping seems
  possible. Thus, it is not just the fact that the diagonal entries of
  $\bm D$ in a Schur decomposition may appear in any order and,
  therefore, $\bm D_1$ and $\bm D_2$ are not well defined in Theorem
  1, Part 1. In turn, the core approach in
  \cite{siavash2021asymmetric}, which relies on a Schur decomposition,
  seems to help little towards solving the targeted problem.
 }

 \textbf{Theorem 1, Part 2.} This result is also incorrect, as the
 proof is flawed at multiple points. First, matrix $\bm Q$ is assumed
 upper triangular, which entails loss of generality given that $\bm Q$
 is upper \emph{quasi-}triangular; see sentence after (5). This
 implies that (8) is not correct and invalidates the rest of the
 proof.


Second, note that the proof relies on the assumption that $\bm T''$ is
full-row rank. Since it has $N_2$ rows, the assumption is $\rank(\bm
T'')=N_2$. However, this violates the Cayley-Hamilton Theorem. To see
this, recall that this theorem states that every matrix $\bm
Z\in \rfield^{N\times N}$ is a root of its own characteristic
polynomial. Let this polynomial be $p(\bm Z)=\sum_{i=0}^N \beta_i \bm Z^i$,
where, as we know, $\beta_N=1$. Then, $p(\bm Z)=\bm 0$ clearly implies that $\bm Z^N =
-\sum_{i=0}^{N-1} \beta_i \bm Z^i$. Multiplying both sides by $\bm
Z^k$, for any integer $k\geq 0$, shows that every power $\bm Z^{N+k}$
is a linear combination of the powers of $\bm Z$ with lower exponents.
This implies that the matrix $[\vect[{\bm Z^0}],\vect[{\bm
Z^1}],\ldots, \vect[{\bm Z^M}] ]$, for an arbitrary $M$, can never be
of rank greater than $N$.

Matrix $\bm T''$ in this paper is obtained by removing rows from
$\bm T= [\vect[{\bm Z^0}],\vect[{\bm Z^1}],\ldots, \vect[{\bm
  Z^{L-1}}] ]$, where $\bm Z=\bm D+\bm Q$. Since $\bm Z$ is
$N \times N$, it follows from the Cayley-Hamilton theorem that
$\rank(\bm T'')\leq \rank(\bm T)\leq N$. However, the proof assumes
that its rank is $N_2=(N^2+N)/2 -m_1-m_2$, which satisfies
$N_2\geq (N^2+N)/2 -N= (N^2-N)/2$. Thus, except for the trivial case
$N\leq 3$, this number is greater than the maximum imposed by the
Cayley-Hamilton Theorem. This contradiction doubly voids the proof.
\changed{Besides, the sufficient condition in Theorem 1, Part 2, is therefore
  useless since it will never hold for $N\geq 3$.}


 \textbf{Centralized algorithm.} The matrix $\bm S$ solving (11) and (12) is
claimed to provide an exact subspace projection; cf. Sec. II and
discussions around (11) and (12). However, such a claim is clearly
false. The cause of this issue is that the condition provided by
Theorem 1, Part 1, is taken as a \emph{sufficient} condition, whereas what
that theorem actually states is that this is a \emph{necessary}
condition. One could however think that the solution actually yields
the desired subspace projection despite the fact that the
argumentation is incorrect. Unfortunately, this is not the case. To
see this, note that given the choices of $\bm W_\parallel$, $\bm
W_\perp$, and $\bm W$ in the paper, an exact projection is attained
only if (7) holds. However, the solution to (11) and (12) will not
generally satisfy (7). To see this, note that (7) comprises $N^2$
equations and these equations are not generally satisfied by an
arbitrarily selected matrix.
In other words, \emph{one would need to somehow constrain $\bm S$ in
these problems to satisfy (7)}. However, this is not the case.


Finally, although not mentioned in \cite{siavash2021asymmetric}, the
technique to deal with  constraint (11f) therein was proposed
in \cite{romero2021fast}.

\textbf{Decentralized algorithm.} In  Equation (13), $\bm \Psi$ is not defined. 
After (13), matrix $\bm \Delta_n$ is improperly defined and
$\epsilon_1$ is not defined. This and other undefined symbols limits
the reader's ability to understand and verify the contents of
Sec.~V.

After some guessing regarding the symbol definitions, it is easy to
see that the algorithm proposed in Sec.~V \changed{need not converge}
to the solution of the targeted problem. To this end, note that
$\bm a_n$ in (13d) is a function of $\bm D^{(l)}$ with $l\neq n$,
where $\bm D^{(l)}$ is an optimization variable of the $l$-th node.
However, $\bm a_n$ is treated as a constant by the algorithm: this
would need to be fixed by introducing $N-1$ constraints similar to
(14e) but for $\bm D^{(l)}$, $l\neq n$, and by introducing $N-1$
additional update equations in (28).  \changed{This clearly means that
  the proposed algorithm is not guaranteed to converge to the solution
  of (12) and, in fact, it is not even clear whether the algorithm
  will converge at all. }

 The fact that this algorithm is wrong is further corroborated by the
simulations, where \emph{the centralized and decentralized algorithms
yield different results}, something that would not be possible if the
algorithm was correctly derived and implemented, given that no
approximation is introduced.

Finally, this algorithm is not truly decentralized since all the nodes
need to know the subspace of interest (cf. $\bm W$ shows up in (28))
and the topology of the graph (cf. $\bm \Psi$ shows up in (28)). Thus,
all nodes have global information and can solve the problem on their
own without interacting with others. One could argue that interaction
with other nodes could decrease the computational
complexity. Unfortunately, this is not the case: such a decentralized
approach would have the same complexity order (yet probably
substantially \changed{more} operations in absolute terms) as its centralized
counterpart, cf. Sec. VIII-D, and it would introduce additional
overhead in terms of time, communication, and power consumption.

 \textbf{Theorem 2.} Although not mentioned, this theorem is based on
 Theorem 2 and Corollary 1 in~\cite{romero2021fast}; see also the
 explanation afterwards. In fact, Appendix C
 in \cite{siavash2021asymmetric} is an adaptation of \cite[Appendix
 C]{romero2021fast}, including the same steps and equations up to
 small modifications.
Besides, the proof in \cite{siavash2021asymmetric} is flawed. To see
that, note that (31) is wrong since the matrices in the set
$\mathcal{S}_\text{PF}$ are not necessarily polynomial shift matrices
as per the definition; cf Sec.~IV. A simple counterexample is the zero
matrix, which is in $\mathcal{S}_\text{PF}$ but is not a ``polynomial
shift''. This clearly voids the proof and Theorem 2.

 \textbf{Equation (21c).} Since $\bm d$ is undefined, the problem (21)
cannot be understood.

 \textbf{Proposition 1.} This is also incorrect. The proof suffers
 from the same issues as the proof of Theorem~1; see above.

\textbf{Simulation study.}  It is remarkable that the simulation
 experiments do not compare with the algorithms proposed
in \cite{romero2021fast}, given that the contribution is based on the
benefits of allowing for asymmetric shift matrices. This could have
been easily done by running an experiment with directed graphs (upon
applying \cite{romero2021fast} on a symmetric adjacency matrix that
results from removing the necessary edges) and another experiment with
undirected graphs, thereby assessing the benefits of the proposed
scheme relative to the state of the art. Instead, the proposed methods
are compared against algorithms not designed for subspace projection.

\section{Conclusions}

This short correspondence has disproved Theorem 1, Theorem 2, and
Proposition 1 in \cite{siavash2021asymmetric} by means of
counterexamples and counterproofs. The proposed centralized algorithm,
which aims at solving (12) was shown to be flawed. Similarly, the
proposed decentralized algorithm is also incorrectly derived.
Besides, it was pointed out which core ideas were obtained from
\cite{romero2021fast} without acknowledgment.  Whenever possible, a
correction was proposed. However, for the most part, \changed{the
  adopted approach based on a Schur decomposition seems intractable}
and amendments are most likely impossible.


\printmybibliography

\clearpage

\end{document}